\documentclass[graybox]{svmult}

\usepackage{custdefns}
\usepackage{type1cm}  
\usepackage{makeidx}  
\usepackage{graphicx} 
\usepackage{multicol}  
\usepackage[bottom]{footmisc}

\usepackage{newtxtext}     
\usepackage[varvw]{newtxmath}    
\makeindex           

\begin{document}

\title*{Search for $\Z/2$ eigenfunctions on the sphere using machine learning}

\author{Andriy Haydys\orcidID{0000-0002-8968-6779} and\\ Willem Adriaan Salm\orcidID{0000-0002-9667-1540}}
\institute{Andriy Haydys, Universit\'e Libre de Bruxelles, \email{andriy.haydys@ulb.be}
\and Willem Adriaan Salm \at Universit\'e Libre de Bruxelles, 
\email{willem.adriaan.salm@ulb.be}}

\maketitle

\abstract*{We use machine learning to search for examples of $\Z/2$ eigenfunctions on the 2-sphere.
For this we created a multivalued version of a feedforward deep neural network, and we implemented it using the JAX library. We found $\Z/2$ eigenfunctions for three cases:
In the first two cases we fixed the branch points at the vertices of a tetrahedron and at a cube respectively.
In a third case, we allowed the AI to move the branch points around and, in the end, it positioned the branch points at the vertices of a squashed tetrahedron.}

\abstract{We use machine learning to search for examples of $\Z/2$ eigenfunctions on the 2-sphere.
For this we created a multivalued version of a feedforward deep neural network, and we implemented it using the JAX library. We found $\Z/2$ eigenfunctions for three cases:
In the first two cases we fixed the branch points at the vertices of a tetrahedron and at a cube respectively.
In a third case, we allowed the AI to move the branch points around and, in the end, it positioned the branch points at the vertices of a squashed tetrahedron.}

\section{Introduction}

In this contribution we describe certain multivalued eigenfunctions of the Laplacian on the sphere $S^2$. These objects appear in the blow-up analysis of multivalued harmonic functions, which is a classical topic in the theory of minimal surfaces and go back to Almgren~\cite{Alm83}. 
More recently, Taubes~\cite{Taubes13PSL2C,Taubes13cASD,Taubes16SW4,Taubes17VW,Taubes18KW} proved that multivalued harmonic functions (and their close cousins multivalued harmonic 1-forms and spinors) provide models for diverging sequences in a whole range of gauge-theoretic problems; see also~\cite{HaydysWalpuski15}. 
They are infinitesimal models for branched multiple covers of special Lagrangian submanifolds~\cite{He23} and play a role in certain fibrations of $\rG_2$-manifolds~\cite{Don17}.
Also, they appeared in the so-called thin obstacle problem~\cite{FranchSavin25}, which has its origins in the theory of elasticity. 
All this provides ample evidence that multivalued harmonic functions and hence, multivalued eigenfunctions of the Laplacian on spheres, are interesting objects to study.

We apply the methods of artificial intelligence (AI), especially deep learning, to search for examples of multivalued eigenfunctions.  Although deep learning is a well-known tool in numerical approach to partial differential equations, to the authors' knowledge this is the first time it is applied to multivalued harmonic functions.

To explain the problem better and the place of AI in our approach, a more precise description of multivalued eigenfunctions will be useful. 
We focus below on 2-valued eigenfunctions, which are also known as $\Z/2$ eigenfunctions, and describe this case in some details.  

Thus, pick an even number of distinct points on the sphere $S^2$, say $Z=\{ p_1,\dots, p_{2n}\}$. Consider the real Euclidean line bundle $\cI=\cI_Z\to S^2\setminus Z$, whose monodromy equals $-1$ along any simple loop encircling any of $p_j$ once. 
If $U\subset S^2\setminus Z$ is any simply connected open subset, we can find a trivialization of $\cI|_U$, i.e., a section $s_0$ of $\cI|_U$ of pointwise length $1$. 
Hence, over $U$ any $u\in \Gamma(S^2\setminus Z;\, \cI)$ can be written as $u=u_0 s_0$ so that we can define $\Delta\, u:= \big (\Delta\, u_0\big ) s_0$, where $\Delta$ on the right-hand side denotes the Laplace--Beltrami operator on $S^2$.  It is easy to check that this twisted Laplacian depends neither on $U$ nor $s_0$, thus is globally well-defined. The equation $\Delta\, u=\lambda u$ then makes sense, however, this is not very interesting without any ``boundary condition'' at $p_j$ since the problem is underdetermined. 

It turns out that appropriate boundary conditions can be posed in the following manner. 
We say that a pair $(u, Z)$ is a $\Z/2$ eigenfunction of the Laplacian, if 
\begin{align}
	&\Delta_{S^2} u =\lambda u,\label{Eq_Eigenfn}\\
	& |u|\le Cr_j^{3/2} \text{ and } |du|\le Cr_j^{1/2} \text{ for all } j\in \{1,\dots, 2n\}, \label{Eq_EigenfnBdryConds}
\end{align} 
where $r_j$ is the distance to $p_j$. 
We would like to point out that along with $u$, \emph{the branching set} $Z$ is a variable. This makes the problem both interesting and challenging.

\begin{remark}
	While the notion of $\Z/2$ eigenfunction can be generalized to any dimension, we focus on the 2-sphere here because it is relatively simple and of interest in applications. 
\end{remark}

Some examples of $\Z/2$ eigenfunctions can be found in~\cite{TaubesWu20Examples,He25Existence,FranchSavin25}. We would like to point out that although the existence $\Z/2$ eigenfunctions has been established, neither eigenfunctions themselves nor corresponding eigenvalues are known explicitly except in the simplest case $n=1$.
The main goal of this contribution is to shed some light on properties of $\Z/2$ eigenfunctions and the spectrum of~\eqref{Eq_Eigenfn} and \eqref{Eq_EigenfnBdryConds}.
To this end, we employed deep learning to search for approximate $\Z/2$ eigenfunctions. 
In the case we know that a configuration $Z$ of points on the sphere can be realized as a branching set of a $\Z/2$ harmonic eigenfunction, we searched for the eigenvalues and $\Z/2$ eigenfunction with $Z$ being fixed. 
This way we obtained (approximations of) $\Z/2$ eigenfunctions with branch configurations consisting of vertices of a regular tetrahedron and a cube;
see Sections~\ref{Sec_Tetrahedron}--\ref{Sec_Cube} for details.  
Letting $Z$ vary, we obtained an example, which is described in Section~\ref{Sec_SquashedTetra}.

The source code for the experiments in this article is available via the following link:  \url{https://github.com/wasalm/2-valued-neural-networks}

\begin{acknowledgement}
This project was partially supported by the Universit\'e Libre de Bruxelles
via the ARC grant “Transversality and reducible solutions in the Seiberg--Witten
theory with multiple spinors”. Also, computational resources have been provided by the Consortium des \'Equipements de Calcul Intensif (C\'ECI), funded by the Fonds de la Recherche Scientifique de Belgique (F.R.S.-FNRS) under Grant No. 2.5020.11 and by the Walloon Region.
\end{acknowledgement}

\section{Local description of 2-valued functions}
Out of the box, deep neural networks cannot be used to describe multivalued functions. 
The goal of this section is to define a 2-valued version of these networks that can be implemented using existing machine learning libraries. First, we give two different descriptions of 2-valued functions.

\subsection{Branched coverings}
As often done in geometry one can describe multivalued functions locally using charts. For this we use the observation that for each real line bundle $\mathcal{I}_Z$, there exists a branched double cover $\pi\colon \hat S\to S^2$ that trivializes $\cI_Z$. Then 2-valued functions, i.e., sections of $\cI_Z$, pull back to ordinary functions on $\hat S$  which are equivariant with respect to the deck transformations. Conversely, any equivariant function on $\hat S$ yields a 2-valued function on $S^2$. 

For simplicity, we work with a finite number of charts $\{U_i\}$, where each $U_i$ is simply connected and each $U_i$ can contain at most one branch point. This way the branched double cover over $U_i$ has an explicit description: If $U_i$ does not contain a branch point,
$\pi^{-1}(U_i)$ is diffeomorphic to $U_i\times\{ \pm 1\}$ and the deck transformation acts by $(z, \e)\mapsto (z, -\e)$. Thus, a 2-valued function on $U_i$ can be identified with an ordinary function $U_i\to \R$.

Suppose a chart $U_i$ contains a branch point. We can choose an identification of $U_i$ with an open subset of $\C$ such that the branching point is the origin and the deck transformation acts by $z\mapsto -z$. Thus, 2-valued functions on $U_i$ can be viewed as odd functions on $\C$ defined on open subsets. 

Finally, we can multiply any locally defined 2-valued function with a smooth bump function to make it global. Conversely, any  2-valued function $u$ is a linear combination of a set of 2-valued functions $u_i$ such that each $u_i$ is supported on the open chart $U_i$. Hence, using our understanding of compactly supported 2-valued functions on each open chart, we can create all possible 2-valued functions on $S^2$.

\subsection{Branch cuts}
The above description using branched double covers has the advantage that both the cover map and the deck transformation are smooth. 
So, by using this formalism we can easily create 2-valued functions that are smooth on $S^2 \setminus Z$. Moreover, this formalism has the advantage as we can turn any map on the branched double cover to a 2-valued map by the projection
\begin{equation}
	\label{eq:projection-to-multivalued-function}
	f(x) \mapsto f(x) - f(-x).
\end{equation}

One disadvantage of this method is that does not give a canonical choice for a fundamental domain. That is, there is no prevalent choice for ${z}$ over $- {z}$. If the choice is not consistent between different charts, the resulting function will not be smooth. A second disadvantage of this method is visualization. Namely, the genus of the branched double cover will depend on the number of branch points and the cover map is not easy to understand for the generic case. 

To fix these issues we give an alternative description of 2-valued functions using branch cuts. Namely, consider embedded non-intersecting curves on $S^2$ between the branch points $p_i$ and $p_1$ for all $i \in \{2 \ldots, 2n\}$. The union of all these curves is a star-shaped graph $G$ on $S^2$. The complement of $G$ in $S^2$ is a simply connected open set and it does not contain any branch point. Hence, a 2-valued function on $S^2 \setminus G$ is represented by a smooth function $f$ from $S^2 \setminus G$ to $\R$. Moreover, there is an explicit description of the monodromy of $f$. Namely, when one crosses the branch cut $G$, $f$ must flip sign.

We will use branch cuts for visualization and for a consistent choice of fundamental domain. We still need to use the first description using charts though. Namely, this method has the disadvantage that the monodromy condition is described in terms of limits. Since computers have a finite resolution, calculating limits will often lead to numerical instability.

\subsection{Multivalued feedforward neural networks}
Recall that for a set of natural numbers $n_1, \ldots, n_N$, a \textit{feedforward neural network with $N$ layers} is a function $NN\colon\R^{n_1} \to \R^{n_N}$ of the form
\begin{equation}
	\label{eq:example-deep-neural-network}
	NN =  \sigma_{N} \circ W_{N-1, N} \ldots \circ  W_{2,3}\circ \sigma_2 \circ W_{1,2},
\end{equation}
where $W_{i, \, i+1}\colon \R^{n_i} \to \R^{n_{i+1}}$ is an affine function and $\sigma_i\colon \R^{n_i} \to \R^{n_i}$ is a set of non-polynomial functions called \textit{activation functions}. The affine functions $W_{i,\, i+1}$ are called \textit{parameters} of a neural network and the goal of machine learning is to find a suitable choice of $W_{i,\, i+1}$ for a fixed set of $n_i$ and $\sigma_i$.

\begin{figure}[b]
	\sidecaption
	\includegraphics[width=7.6cm,trim={4cm 0 18cm 4.1cm},clip]{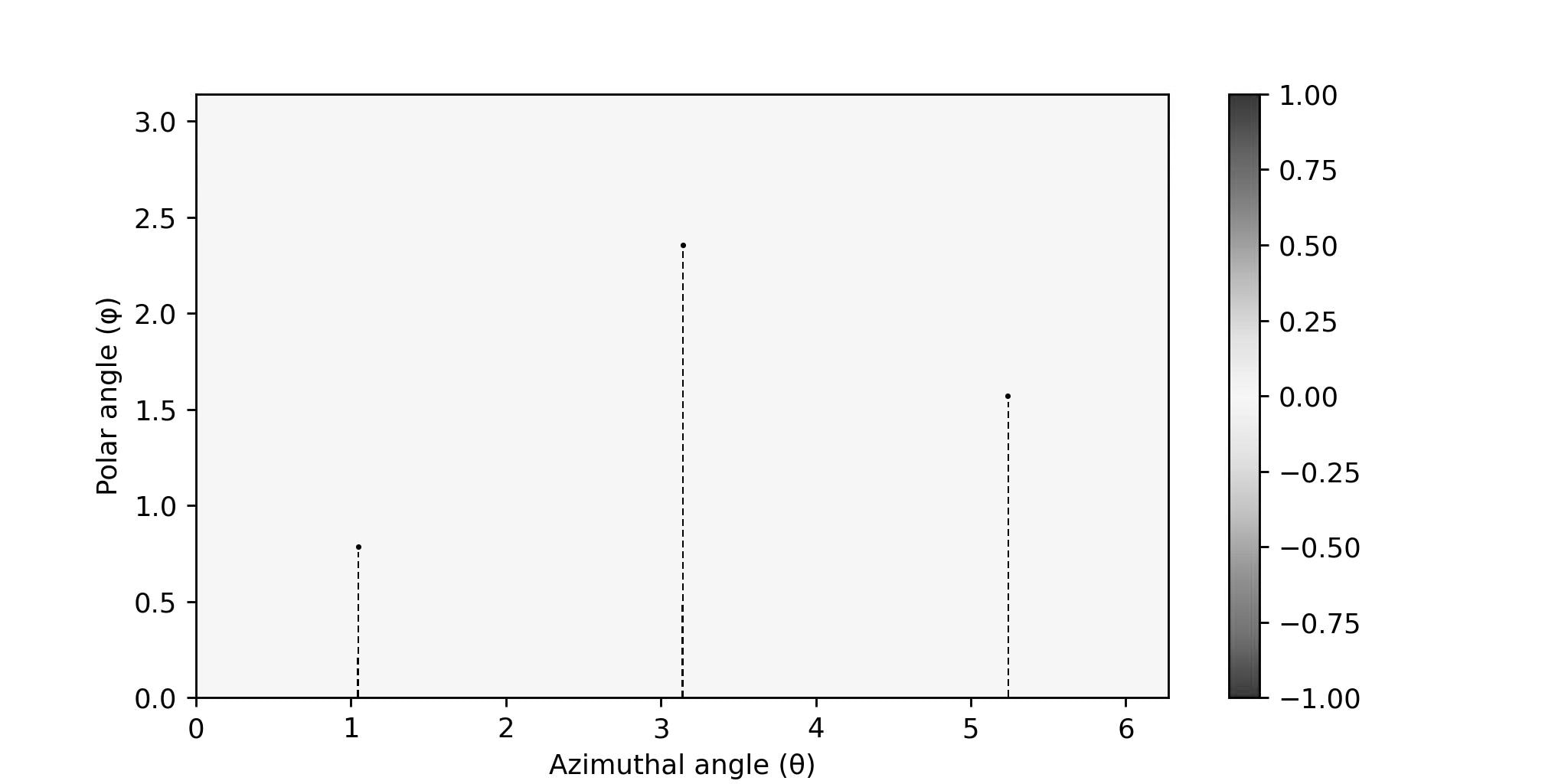}
	\caption{Example of a star-shaped graph $G$ on $S^2$ when displayed using spherical coordinates. The centre of the star is at the north pole ($\phi = 0$) and the other branch points are connected using sections of great circles. In spherical coordinates, these edges correspond to line sections with a constant value of $\theta$.}
	\label{fig:graph}
\end{figure}

With our previous understanding of 2-valued functions, we define our 2-valued feedforward neural network as follows. 
First fix
\begin{enumerate}
	\item a finite set of simply connected charts $U_i$ that cover the 2-sphere, and each chart contains at most one branch point,
	\item a set of smooth bump functions $\chi_i$ such that each $\chi_i$ is supported in $U_i$, and
	\item a neural network $NN_i$ as described in~\eqref{eq:example-deep-neural-network} for each chart $U_i$.
\end{enumerate}
Then for each $x \in S^2$, consider the charts $U_i$ such that $x \in U_i$. Lift $x$ to $z_i$ in the branched double cover of $U_i$ in $\C$. 
If $U_i$ does not contain a branch point, define
$$
u_i = NN_i(z_i) \cdot \chi_i(z_i).
$$
Multiplication with $\chi_i$ ensures that our result is globally defined.
If $U_i$ does contain a branch point, define
\begin{equation}
	\label{eq:formula-2-valued-function-on-chart}
	u_i = (NN_i(z_i) - NN_i(-z_i)) \cdot \chi_i(z_i) \cdot |z_i|^2.
\end{equation}
Multiplication with $|z_i|^2$ ensures that the decay rate to the branch point is of order $\frac{3}{2}$.
Finally, we define a 2-valued feedforward neural network $u$ as the sum of all $u_i$.
\begin{remark}
	In our initial implementation, the choice of fundamental domain was not consistent between different charts. Instead of considering a different lift of $x$, we multiplied the output of $u_i$ with a factor of $\pm 1$. From~\eqref{eq:formula-2-valued-function-on-chart} one can see that this approach is equivalent to our initial one.
\end{remark}
To implement this in a computer program, we had to systematically choose our charts and set of bump functions. For this we assumed that one of our branch points was at the north pole, and we picked this as the centre of our star-shaped graph $G$. For the edges of the graph $G$, we used arcs of great circles. For this we had to exclude the degenerate cases where these arcs overlap. An example of a non-degenerate graph is shown in Figure \ref{fig:graph}.

\begin{figure}[b]
	\sidecaption
	\includegraphics[width=7.6cm,trim={4cm 0 18cm 4.1cm},clip]{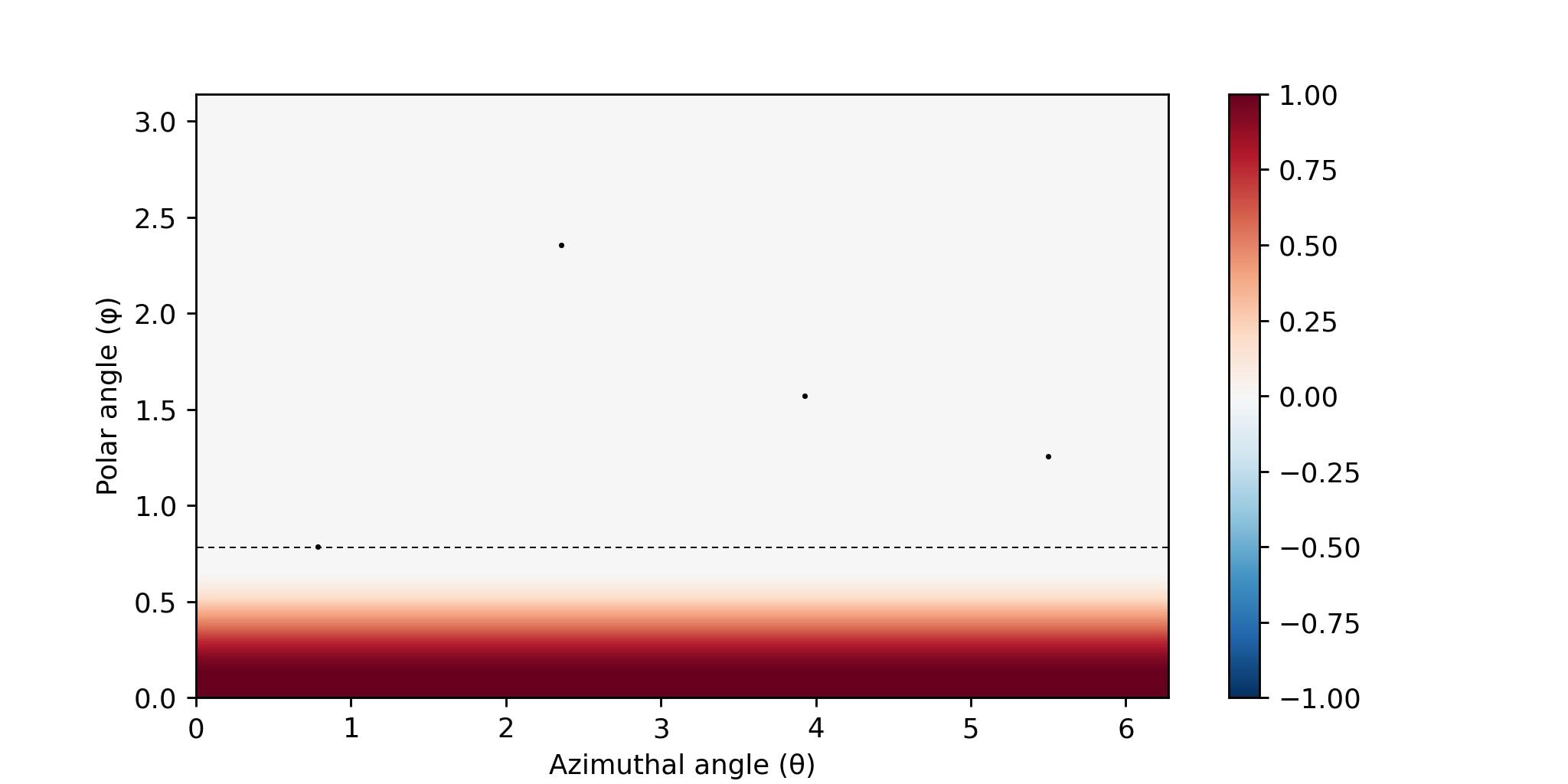}
	\caption{The chart centred at the north pole $\phi = 0$. In this picture, this is the region below the dotted line. On $S^2$ this chart can be visualised as a spherical dome centred at the north pole.
	On this chart we also defined a bump function, and the value of this function is depicted by the red gradient scale. The remaining dots are the locations of the generic branch points.
	}
	\label{fig:chart-north-pole}
\end{figure}

To systematically pick charts, we first considered a spherical dome around the north pole. This is shown in Figure \ref{fig:chart-north-pole}. As seen in the picture, we choose our dome as large as possible without adding an extra branch point. To understand the necessity of this, first recall that given a function $f$ and a bump function $\chi$ with $\|\chi\|_{C^0} = 1$, we have
\begin{equation*}
	\Delta (f \chi) 
	= f \Delta \chi + \chi \Delta f - 2 \langle \operatorname{d} f , \operatorname{d} \chi \rangle.
\end{equation*}
If the support of $\chi$ is concentrated in a small region, the $C^2$-norm of $\chi$ will be large.
Hence, in this case $\Delta \chi$ and $\operatorname{d} \chi$ will be large and this will heavily contribute to the error. 
Thus, to minimize this contribution, it is important to choose the support of $\chi$ as large as possible. 

\begin{figure}[b]
	\sidecaption
	\includegraphics[width=7.6cm,trim={4cm 0 18cm 4.1cm},clip]{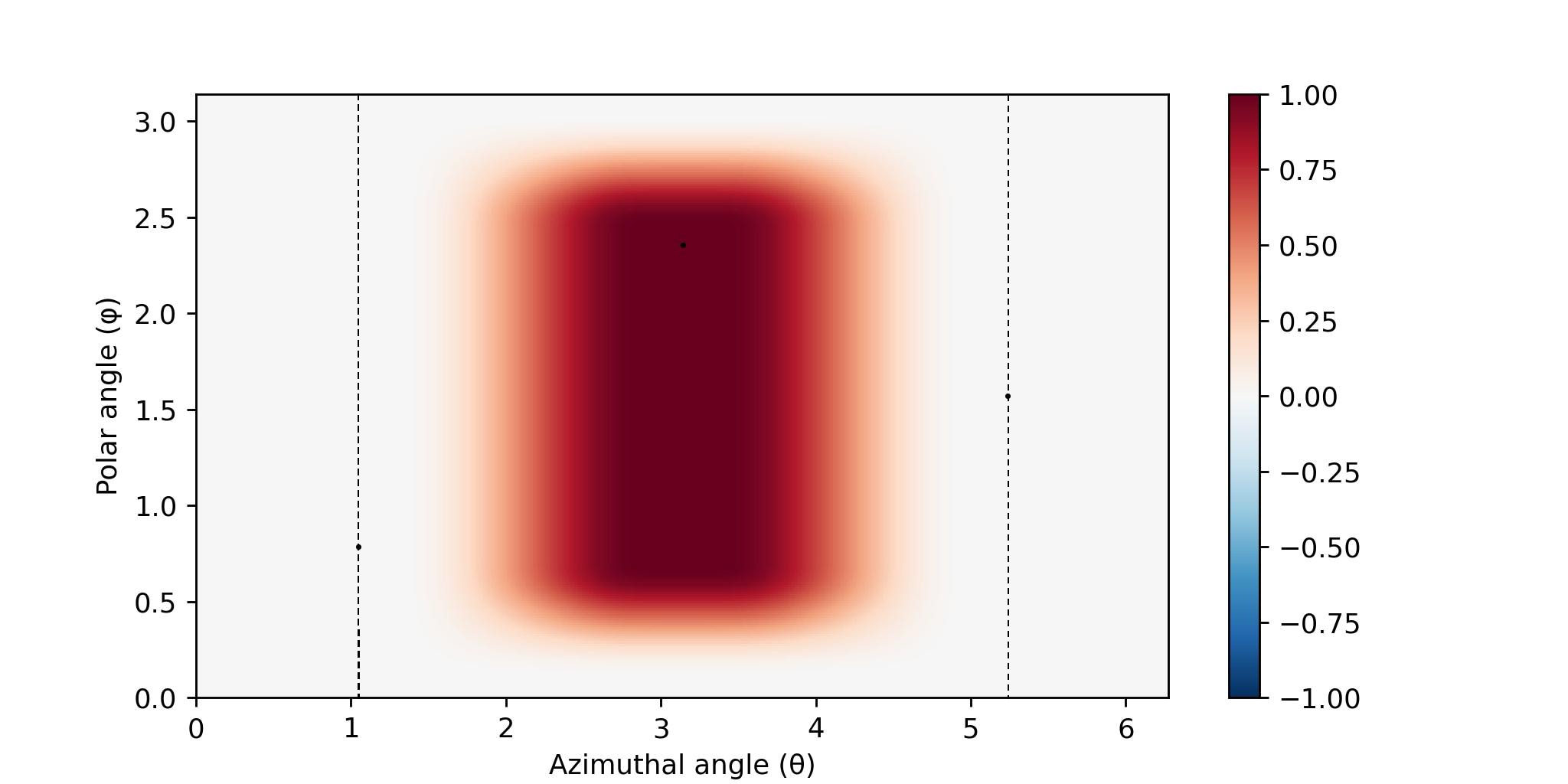}
	\caption{The chart centred at a generic branch point. In this picture, this is the region between the two dotted lines. On $S^2$ this chart has the shape of an orange slice.
	On this chart we also defined a bump function, and the value of this function is depicted by the red gradient scale. The remaining dots are the locations of the generic branch points.}
	\label{fig:chart-generic}
\end{figure}

For a generic branch point, we consider a chart in the shape of an orange slice. This is shown in Figure \ref{fig:chart-generic}. Explicitly, we consider spherical coordinates $(\phi, \theta)$ and we calculated the azimuthal angle $\theta_{\pm}$ for the first branch point east and west of the branch point in question. The domain of our chart is then  $\{\theta \in S^1, \colon \theta_-  < \theta < \theta_+ \}$.

\begin{figure}[b]
	\sidecaption
	\includegraphics[width=7.6cm,trim={4cm 0 18cm 4.1cm},clip]{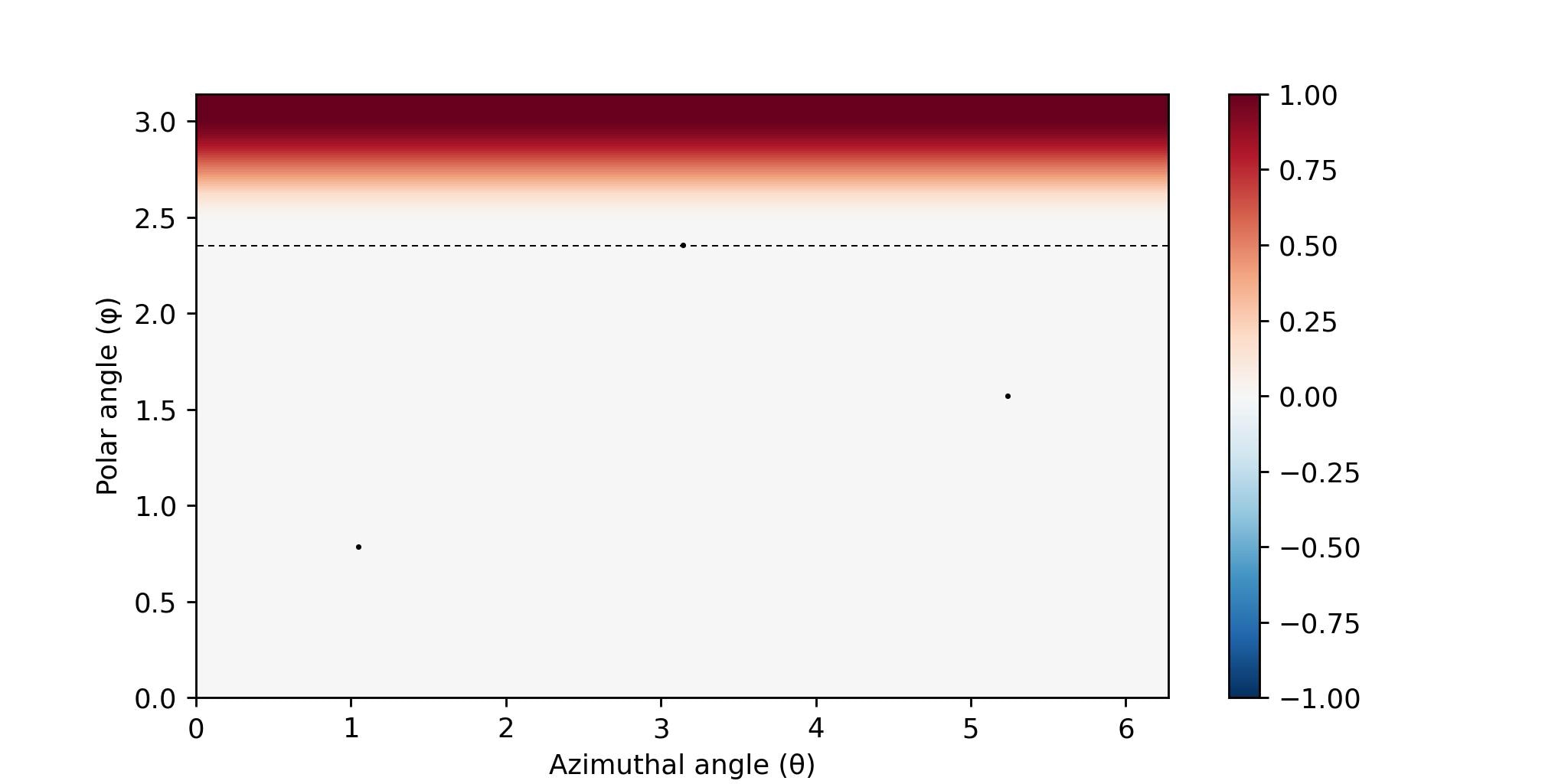}
	\caption{The chart centred at the south pole when the south pole does not contain a branch point. In this picture, this is the region above the dotted line. On $S^2$ this chart can be visualised as a spherical dome centred at the south pole.
	On this chart we also defined a bump function, and the value of this function is depicted by the red gradient scale. The remaining dots are the locations of the generic branch points.}
	\label{fig:chart-south-pole2}
\end{figure}

\begin{figure}[b]
	\sidecaption
	\includegraphics[width=7.6cm,trim={4cm 0 18cm 4.1cm},clip]{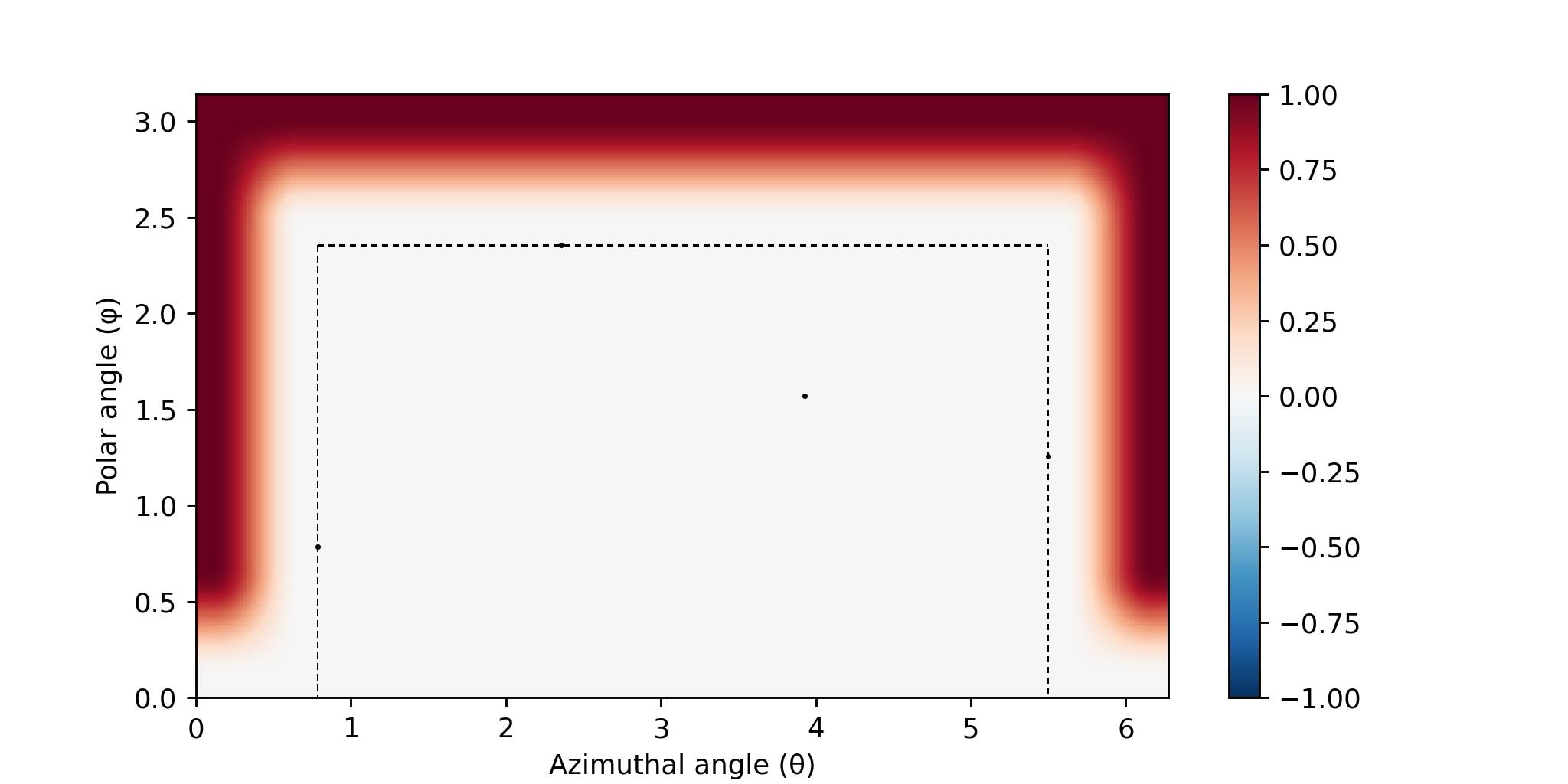}
	\caption{The chart centred at the south pole when the south pole does contain a branch point. In this picture, this is the region outside the dotted rectangle. On $S^2$ this chart can be visualised as the union of a spherical dome and an orange slice.
	On this chart we also defined a bump function, and the value of this function is depicted by the red gradient scale. The remaining dots are the locations of the generic branch points.}
	\label{fig:chart-south-pole1}
\end{figure}

We also add a chart centred around the south pole. In our program we considered two cases depending on whether there is  a branch point at the south pole $\phi = \pi$. If  there is no branch point at the south pole, then we can use a spherical dome like we did with the north pole. This is shown in Figure \ref{fig:chart-south-pole2}. In the case that there is a branch point at the south pole, we assumed that the edge connecting the north and south pole is at $\theta = 0$. To cover the south pole and this edge, we consider a chart that is a union of a spherical dome and an orange slice. This is shown in Figure~\ref{fig:chart-south-pole1}.


Finally, we constructed a set of bump functions that are supported on each chart. For this we picked the standard smooth step function $\mathcal{S}\colon \R \to [0,1]$ that is defined as
\begin{equation}
	\label{eq:smooth-step-function}
	\mathcal{S}(x) = \frac{\mathcal{T}(x)}{\mathcal{T}(x) + \mathcal{T}(1-x)} \, ,\qquad\text{where} \qquad \mathcal{T}(x) = \begin{cases}
	e^{-1/x} & x > 0, \\
	0  & x \le 0.
\end{cases} 
\end{equation}
For each boundary segment in Figures \ref{fig:chart-north-pole}, \ref{fig:chart-generic}, \ref{fig:chart-south-pole2} and \ref{fig:chart-south-pole1} we considered a rescaled version of $\mathcal{S}$. To combine these step functions into a bump function we multiplied them together.
\begin{remark}
	This step function $\mathcal{S}$ was not predefined in our software. 
	To make it smooth at zero, we used a trick described in \cite{JaxFaq}.
\end{remark}

\subsection{The loss function}
With this explicit description of a 2-valued feedforward neural network, we can use standard machine learning techniques to approximate a $\Z_2$ eigenfunction. For a fixed number of branch points, we define a \textit{loss function}
$$
L \colon 
 \left\{ \parbox{5em}{2-valued deep neural networks}
 \right\}
 \times
 \left\{ \parbox{5em}{Possible eigenvalues}
 \right\} \to [0, \infty ) 
$$
such that $(u, \lambda)$ is a minimiser of $L$ if and only if $u$ is a 2-valued eigenfunction of the Laplacian  with eigenvalue $\lambda$ on $S^2$. Using stochastic gradient descend we can search for a minimiser.
To do this, we first considered a small neighbourhood $U$ of the graph $G$ and we defined the loss function $L(u, \lambda)$ as a weighted sum of the following items:
\begin{enumerate}
	\item The $L^2$ error $\|\Delta_{S^2} u - \lambda \: u\|_{L^2(S^2 \setminus U)}$. 
	\item \label{It_C0Error} The weighted $C^0$ error 
	$\|w \cdot(\Delta_{S^2} u - \lambda \: u)\|_{C^0(S^2 \setminus U)}$, where $w$ is a differentiable function on $S^2$ that equals $|z|$ close to the branch points and is constant anywhere else. 
	\item \label{It_Norm1} The normalisation error 
	$(1 - \|u\|_{C^0(S^2 \setminus U)})^2$. 
	\item \label{It_BranchPtsCoin} An error term that penalizes if two branch points are getting too close or 
	when two edges of $G$ almost overlap.
	\item Optionally, the eigenvalue $\lambda$. 
\end{enumerate}

Notice that to approximate the $L^2$ norm of $\Delta_{S^2} u - \lambda \: u$, we picked a random set of points on $S^2 \setminus U$ and calculated the root mean square error. The reason we need to exclude a neighbourhood of the graph $G$ is due to numerical stability. Namely, in training we sometimes received a division-by-zero error while we calculated the Laplacian in this region.

Concerning~\ref{It_C0Error}:  The addition of the weight function $w$ is needed as our 2-valued neural networks decay with order $|z|^{3/2}$ near the branch points. Hence, the Laplacian of our 2-valued neural networks can blow up with the rate $|z|^{-1/2}$ and this weight function is required to compensate for this blow up. To approximate the $C^0$ norm of $w \cdot(\Delta_{S^2} u - \lambda \: u)$, we considered the same random set of points as before and we calculated the maximum absolute error.

Concerning~\ref{It_Norm1}: This is added to normalise $u$, because any scalar multiple of a $\Z_2$ eigenfunction is also a $\Z_2$ eigenfunction. Alternatively, one can divide the above terms by $\|u\|$ and ignore this loss term. However, adding this term helps with the numerical stability of our implementation.

Concerning~\ref{It_BranchPtsCoin}:  If two branch points almost coincide or if two edges of the graph $G$ almost overlap, our implementation of a 2-valued feedforward neural network degenerates and our program gives division-by-zero errors. In our software, we defined a `safe distance' $d_0$ and we measured the minimum distance\footnote{As this distance has no significant meaning, we didn't use the metric on $S^2$, but we just calculated the difference in azimuthal angle and used this as a distance measure instead.} $d$ between the branch points and between the edges. As our penalty, we considered
$$
\left(
	\frac{1}{d} - \frac{1}{d_0}
\right) \cdot \left(1 - \mathcal{S} \left(
\frac{2d - 1}{d_0}
\right)\right)
$$
where $\mathcal{S}$ is the smooth step function defined in Equation \ref{eq:smooth-step-function}. This way the penalty will be of order $\frac{1}{d}$ when the distance $d$ is small and zero of $d > d_0$.

Finally, adding $\lambda$ to the loss function results in searching for the eigenfunction with the smallest eigenvalue.

We would like to point out in passing, that the idea to combine the $L^2$ norm with the $C^0$ norm is not new. According to \cite{Bar19,Bar21}, the addition of the $C^0$ norm controls the outliers that the $L^2$ norm cannot find.

\subsection{Implementation}
With this mathematical background in mind, we implemented a machine learning algorithm in python. Due to the non-standard, mathematical nature of our 2-valued neural networks, we preferred to use the JAX~\cite{JaxHome} library over TensorFlow and PyTorch. This way we were able to describe our problem as if we were using NumPy, while our code is automatically optimized for the use on a GPU. Also, the feature that JAX can differentiate functions was an extra bonus for us.

The downside of this library was that standard ML features, like many standard optimisers, were not implemented. However, in the documentation~\cite{JaxTutorial} there was an example for the implementation of stochastic gradient descend. In the end we had to write around 100 lines of code reimplementing the needed standard features.

In summary, our machine learning algorithm makes use of the following machine learning techniques:

\begin{itemize}
	\item Since we wish to avoid local minima, we used stochastic gradient descent. We ran for 5000 epochs, and we split our dataset of 524,288 random points on $S^2$ in 256 mini-batches. Between rounds we shuffled the dataset over the mini batches. 
	\item To prevent overfitting, we considered a separate set of 8192 points. We configured our 'hyperparameters' using this set. Finally, we evaluated all errors using a third dataset of 16384 points. All errors in this paper are measured using this last dataset.
	\item For each feedforward neural network we assigned 3 hidden layers with 64 neurons. Between the layers we used the GELU function~\cite{GeLU} component-wise as the activation function.
	\item In order to ease to a correct solution, we used a cosine learning rate scheduler and every 1000 rounds we performed a warm restart. Namely, after every 1000 rounds, we set the learning rate to 0.001 and over time we decreased it to 0.0001.
	\item For optimizer we used the AdamW algorithm~\cite{AdamW} with parameters $\beta_1 = 0.9$, $\beta_2 = 0.999$, $\epsilon = 10^{-8}$ and $\lambda = 0.004$.
	\item To initialise the parameters of the neural networks we use the method described in ~\cite{He15}. That is, we initialise the bias\footnote{That is the constant part of an affine function $W$.} to zero, and pick the weights\footnote{That is the linear part of an affine function $W$.} randomly from
	$\mathcal{N}\left(0,\sqrt{\frac{2}{n_{l-1}}}\right)$ where $n_l$ is the number of neurons in the $l$'th layer. The eigenvalue and the position of the branch points are initialised by hand.
\end{itemize}

For hardware we were thankful to use the cluster provided by the Consortium des \'Equipements de Calcul Intensif (C\'ECI), which provided us access to an NVIDIA RTX 6000 ADA.

\section{Results}

\subsection{The tetrahedral case}
\label{Sec_Tetrahedron}

We searched first a 2-valued eigenfunction with the branch points located at the vertices of a regular tetrahedron. We chose this as its existence was shown by Taubes and Wu~\cite{TaubesWu20Examples}. To find this example, we denied the AI the possibility to move the branch points around. After some trial and error, we found some suitable weights for the loss function: We decided to weigh the $L^2$-loss with a factor of 10, the $C^0$-loss with a factor of 1 and the eigenvalue with a factor of 2. To prevent the system to converge to the zero map, we weighted the $C^0$-norm with a factor of 100.

\begin{figure}[b]
	\sidecaption
	\includegraphics[width=7.6cm,trim={1.0cm 0 2.5cm 1cm},clip]{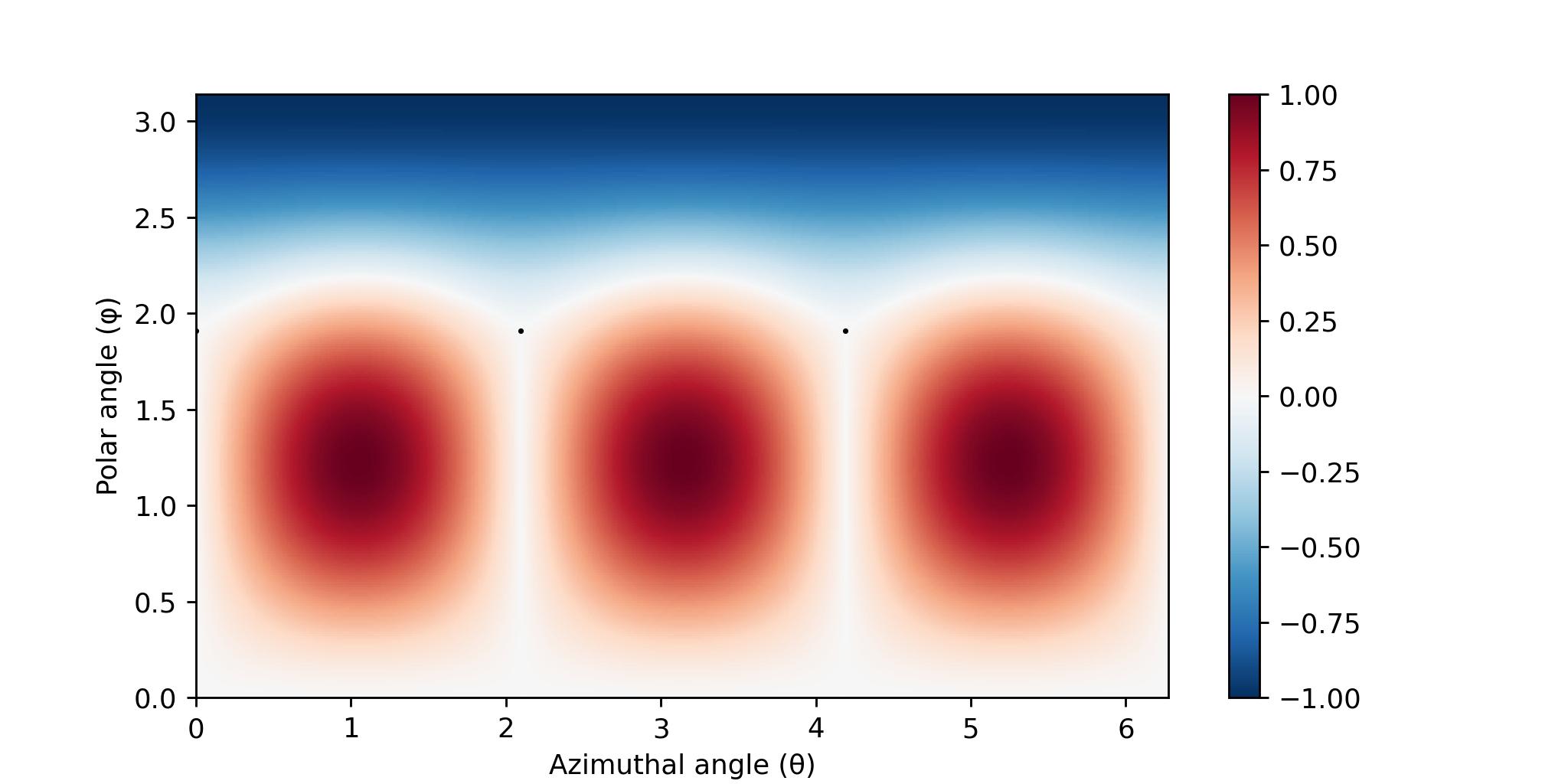}
	\caption{Plot in spherical coordinates of a 2-valued eigenfunction with the branch points located at the vertices of a tetrahedra. One branch point is located at the north pole ($\phi = 0$); Locations of the other branch points are marked by  dots.}
	\label{fig:plot-tetrahedra}
\end{figure}

\begin{figure}[b]
	\includegraphics[width=\textwidth]{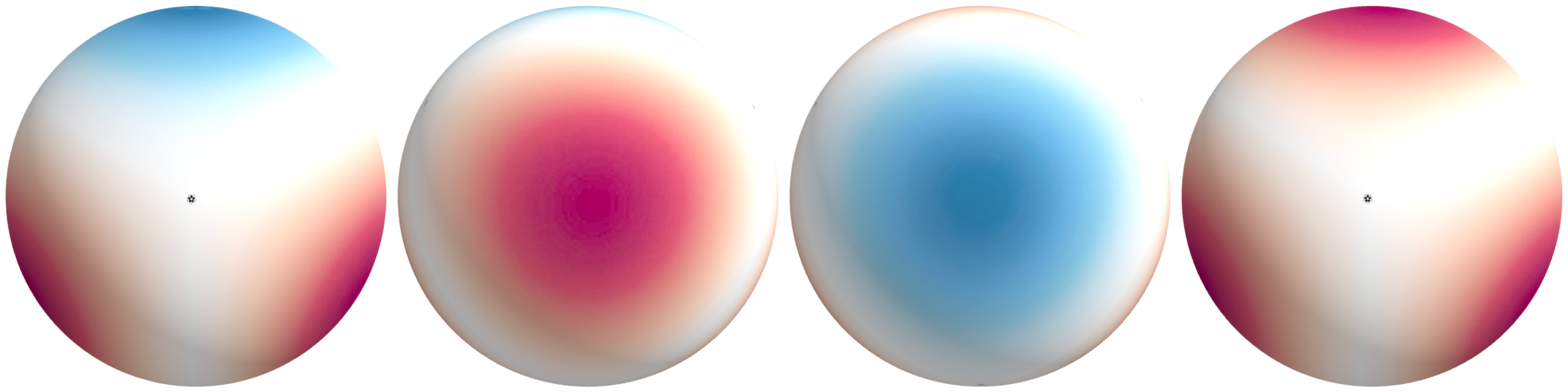}
	\caption{3d renders of a 2-valued eigenfunction with the branch points located at the vertices of a regular tetrahedron.}
	\label{fig:3d-tetrahedra}
\end{figure}

Figure \ref{fig:plot-tetrahedra} shows the 2-valued eigenfunction we have found. It has an expected eigenvalue of $5.154$, the root mean squared error of $0.017$ and  a weighted maximal absolute error of $0.065$. Using a 3d render\footnote{An interactive version is available in our code repository.} shown in Figure \ref{fig:3d-tetrahedra}, we studied its symmetries. Even though we didn't force it, the solution looks invariant under the isometries of a tetrahedron. 
It appears that we have found an approximation for the example, whose existence was established in~\cite{TaubesWu20Examples}. 

Observe that the plot contains the following useful piece of information. 
Notice first that near any branch point $x\in S^2$, a $\Z/2$ eigenfunction admits an asymptotic expansion of the form
$$
	u=\Re  \big ( a \, w^{3/2} \big ) + O\big (|w|^{5/2}\big ),
$$
where $a$ is a complex number and $w$ is a complex coordinate on $S^2$ centred at $x$. 
This implies that $a\neq 0$ if and only if there are exactly three branches of the zero locus of $u$ meeting at $x$. 
If the coefficient $a$ in the above expansion does not vanish at each branch point, the eigenfunction $u$ is called non-degenerate. 
These non-degenerate eigenfunctions play a particular role in the theory of $\Z/2$ harmonic functions, see~\cite{ChenHe24, Don21Deform, HMT23, Takahashi23Moduli} for more details.    

Coming back to the plot, the zero locus of $u$ is represented by white lines and in the example under consideration contains exactly three branches near each branch point.
This is in line with~\cite[Prop.\, 5.17]{ChenHe24}.  

\begin{figure}[b]
	\sidecaption
	\includegraphics[width=7.6cm,trim={1.0cm 0 2.5cm 1cm},clip]{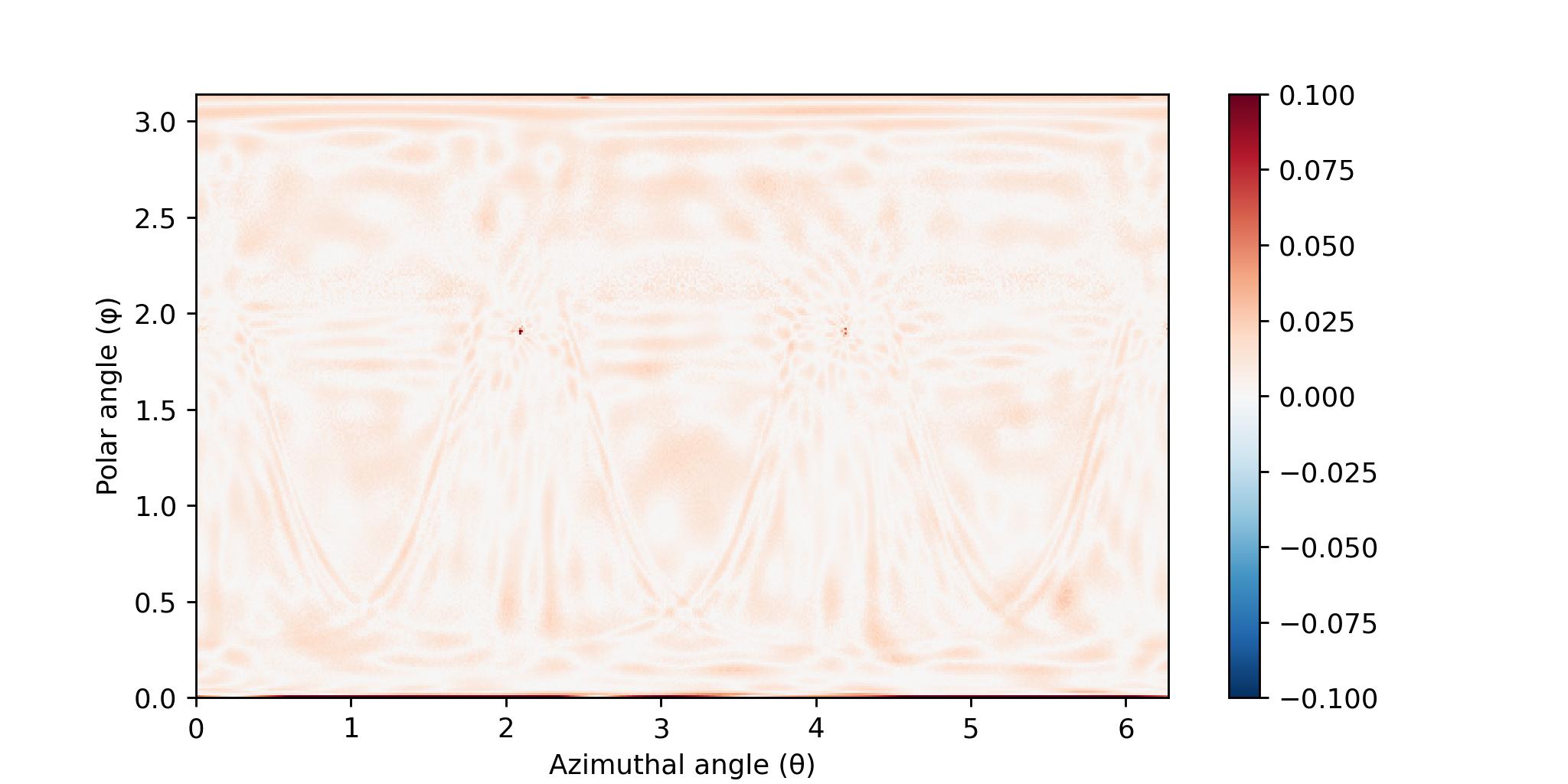}
	\caption{Plot in spherical coordinates of the weighted error for the 2-valued eigenfunction with the branch points located at the vertices of a tetrahedron. 
	Due to numerical instability, aggravated by the fact that our program uses 32-bit floating numbers, the error is concentrated at the branch points.}
	\label{fig:error-tetrahedra}
\end{figure}

A plot of the error in this example is shown in Figure~\ref{fig:error-tetrahedra}. 
As one can see, the error is concentrated in a neighbourhood of the branch points. This is due to numerical instability because we are multiplying a value of order $|z|^{-1/2}$ with that of order $|z|$.
This error plot gave us the confidence that we found a 2-valued eigenfunction even while we exclude small neighbourhoods around our branch points from our testing set. We did the same for the upcoming examples and we were able to make the same conclusion.

\subsection{The cubic case}
\label{Sec_Cube}
Another example we searched with our AI was a 2-valued eigenfunction where the branch points are the vertices of a cube. Again, its existence was proven by Taubes and Wu~\cite{TaubesWu20Examples}. The settings were almost identical to the tetrahedral case; however, the convergence of the eigenvalue was fairly slow. Therefore, we changed the weight of the eigenvalue in the loss function to $0.25$.

\begin{figure}[b]
	\sidecaption
	\includegraphics[width=7.6cm,trim={1.0cm 0 2.5cm 1cm},clip]{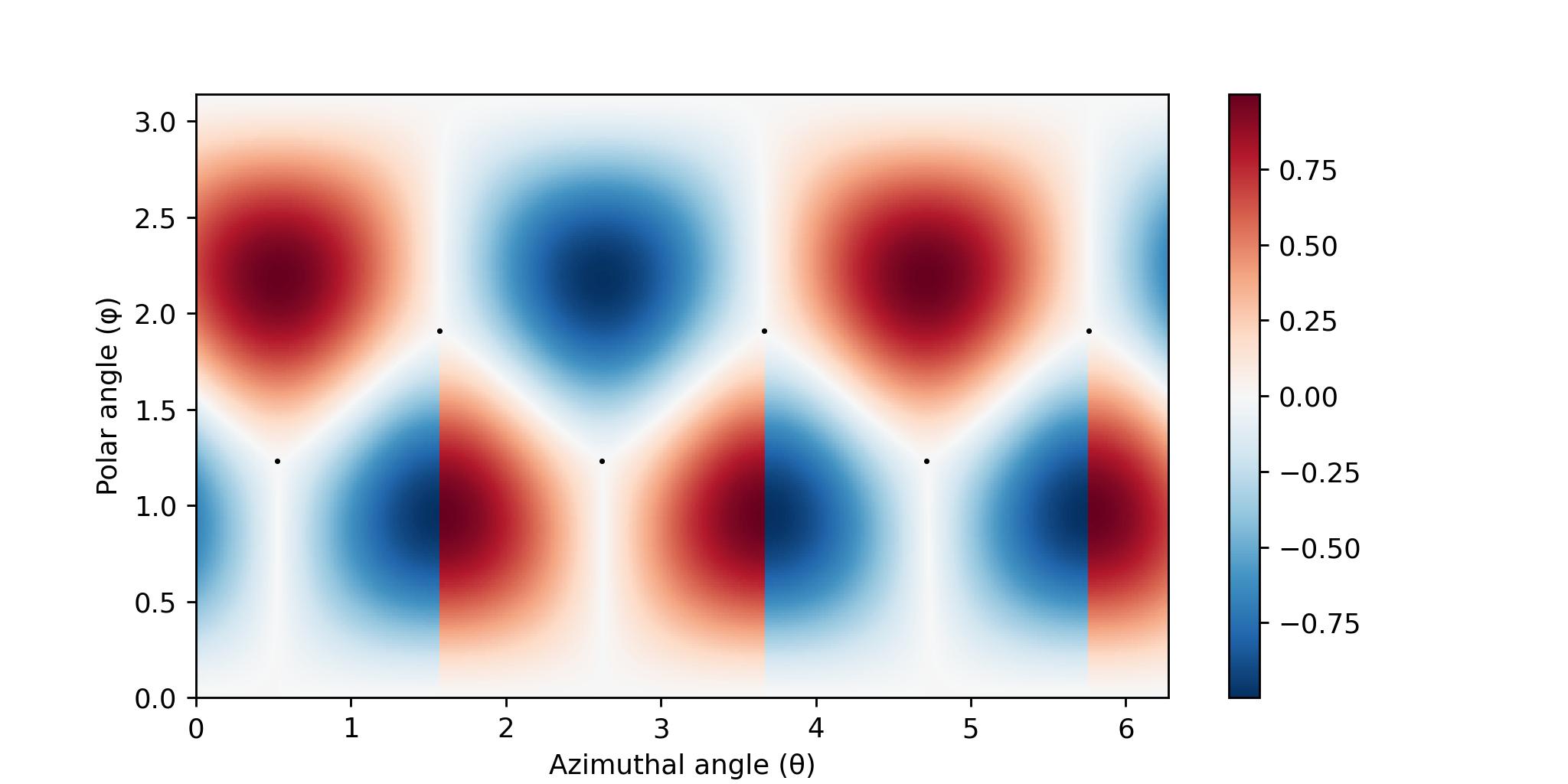}
	\caption{Plot in spherical coordinates of a 2-valued eigenfunction with the branch points located at the vertices of a cube. Two branch points are located at the poles ($\phi \in \{0, \pi\}$) and the other branch points are located at the dots.}
	\label{fig:plot-cube}
\end{figure}

\begin{figure}[b]
	\includegraphics[width=\textwidth]{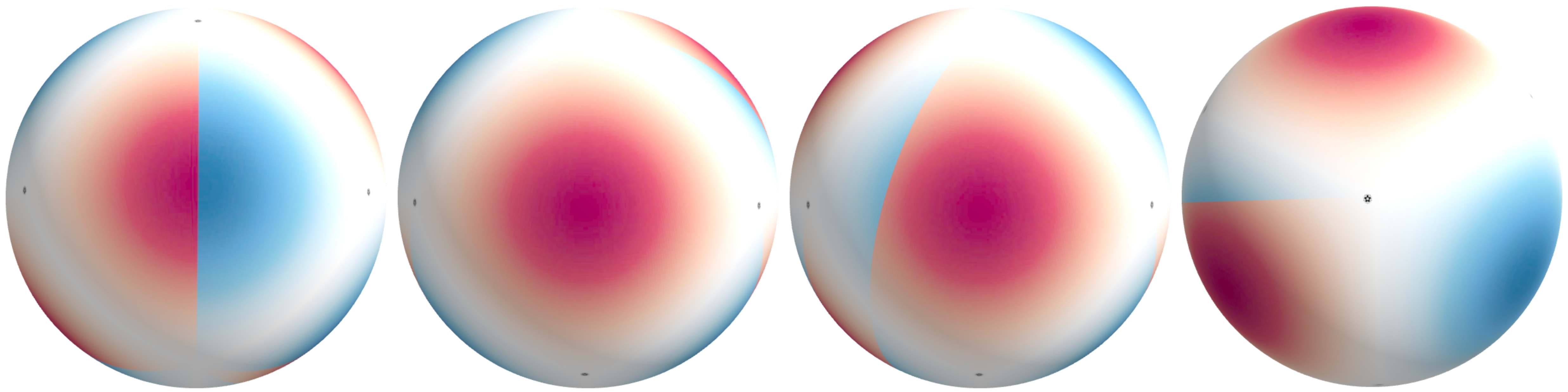}
	\caption{3d renders of a 2-valued eigenfunction with the branch points located at the vertices of a cube.
	}
	\label{fig:3d-cube}
\end{figure}

Figure \ref{fig:plot-cube} shows the 2-valued eigenfunction we have found. It has an expected eigenvalue of $8.098$, the root mean squared error if $0.117$ and a weighted maximal absolute error of $0.462$. 3D renders are shown in Figure \ref{fig:3d-cube}, which again look invariant under the isometries of cube. 

Notice that just as in the tetrahedral case, this example is also non-degenerate.  

\subsection{The squashed tetrahedral case}
\label{Sec_SquashedTetra}

Finally, we searched with our AI for a random example. This time we added position of the branch points as one of our parameters of our network and we let the algorithm search for the best locations.

This turned out to be more challenging than we expected. Sometimes we found solutions with $|z|^{1/2}$ decay. Sometimes, all branch points coalesced at the north-pole and often mitigations resulted with solutions with high error rates.

After implementing the following changes, we found a solution. First, we disabled the minimization of the eigenvalue in the loss function. Although this was useful for comparison with ~\cite{TaubesWu20Examples}, it is unlikely that for a fixed set of branch points the eigenfunction with the smallest eigenvalue has a decay rate of $|z|^{3/2}$. Secondly, we forced the minimal distance between branch points to be large. This way we stayed far away from the local minima we easily got trapped in. Finally, we multiplied this penalty with a factor of 5 in our loss function to make it significant.

\begin{figure}[t]
	\sidecaption
	\includegraphics[width=7.6cm,trim={1.0cm 0 2.5cm 1cm},clip]{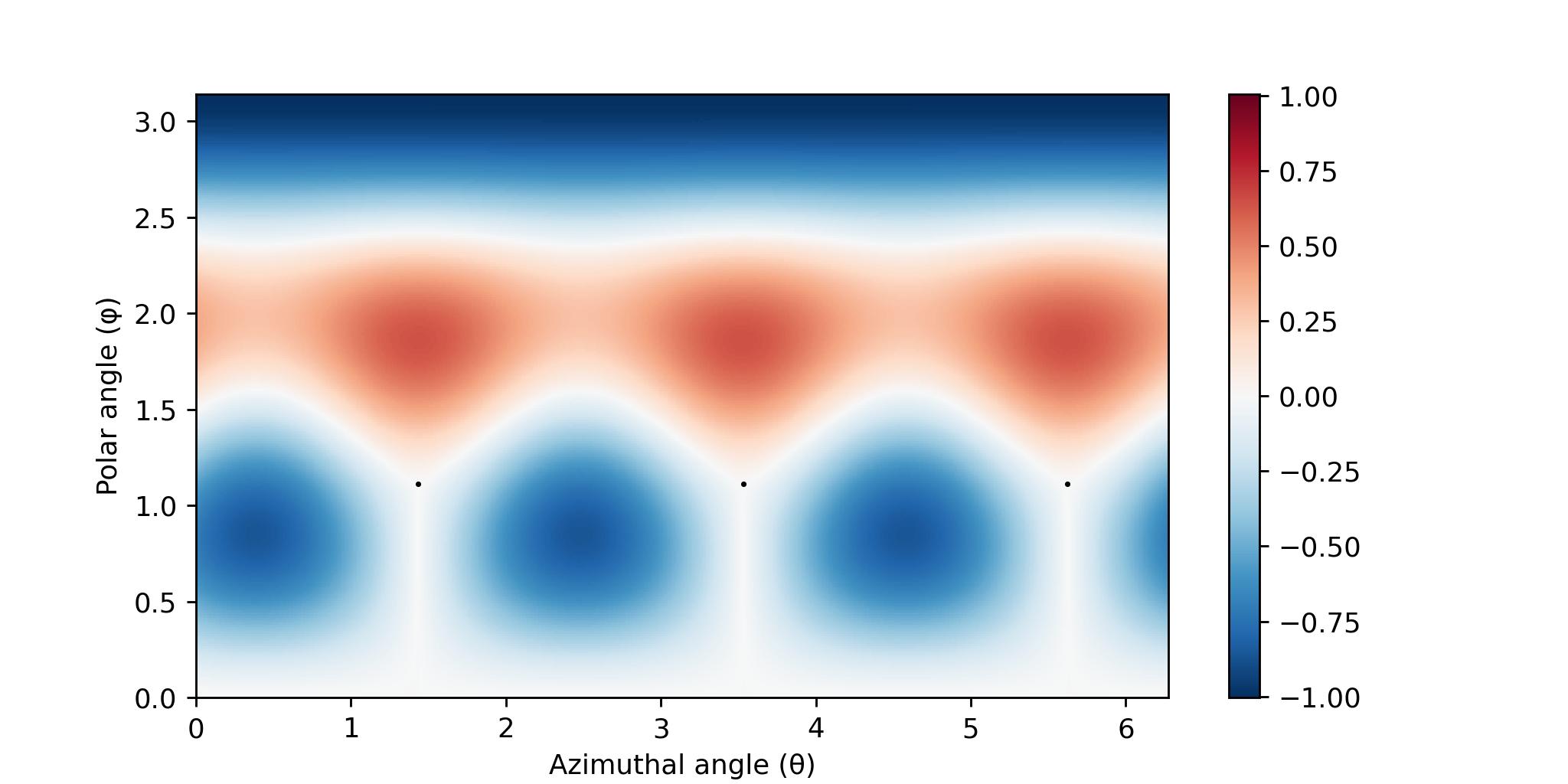}
	\caption{Plot in spherical coordinates of a 2-valued eigenfunction. One branch point is located at the north pole ($\phi = 0$) and the other branch points are located at the dots.}
	\label{fig:plot-new-example}
\end{figure}

Figure \ref{fig:plot-new-example} shows the approximation of the non-degenerate 2-valued eigenfunction we have found. It has an expected eigenvalue of $9.874$ and the root mean squared error is $0.030$ and its weighted maximal absolute error is $0.100$. Notice that albeit there are four branch points, they do not form a \emph{regular} tetrahedron. Up to an error of 0.5\%, the branch points form a regular pyramid with an equilateral triangle as its base. This is only invariant under the dihedral group. The existence of $\Z/2$ eigenfunctions with such branching locus and symmetries was established in~\cite[Sect.\, 6]{ChenHe24}.  We believe our solution approximates one of these 2-valued eigenfunctions.

\begin{figure}[t]
	\includegraphics[width=\textwidth]{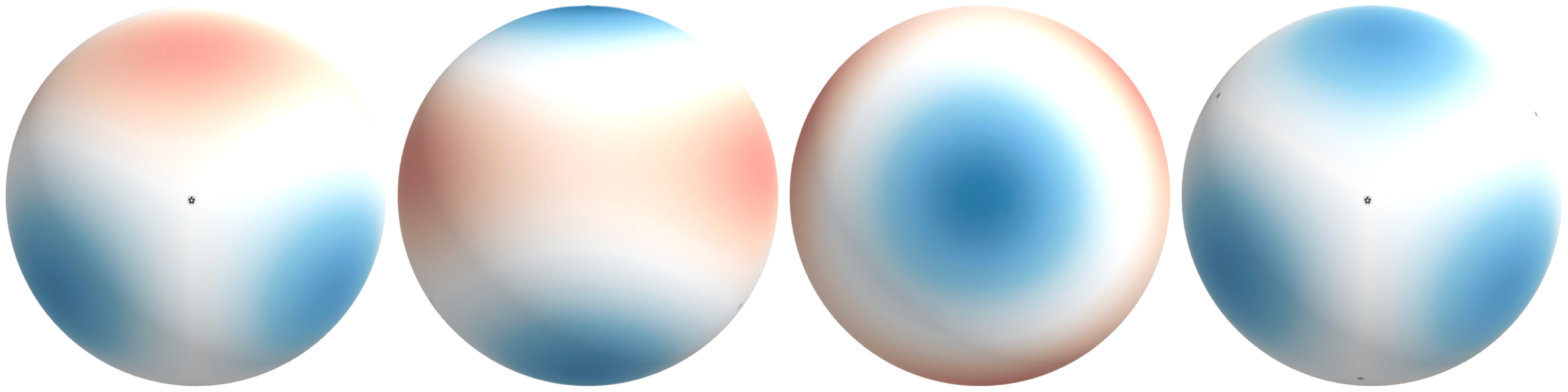}
	\caption{3d renders of the new 2-valued eigenfunction.}
	\label{fig:3d-new-example}
\end{figure}

%
%

\end{document}